\documentclass[12pt]{amsart}
\usepackage{amscd}
\usepackage{amssymb}

\def\Ad{\mathop{\mathrm {Ad}}\nolimits}
\def\ad{\mathop{\mathrm {ad}}\nolimits}
\def\Lie{\mathop{\operatorname {Lie}}\nolimits}
\def\Aff{\mathop{\mathrm {Aff}}\nolimits}
\def\aff{\mathop{\mathrm {aff}}\nolimits}

\def\Id{\mathop{\mathrm {Id}}\nolimits}
\def\Ad{\mathop{\mathrm {Ad}}\nolimits}
\def\ad{\mathop{\mathrm {ad}}\nolimits}

\def\Id{\mathop{\mathrm {Id}}\nolimits}
\def\Aut{\mathop{\mathrm {Aut}}\nolimits}

\def\sgn{\mathop{\mathrm {sgn}}\nolimits}
\def\Lie{\mathop{\mathrm {Lie}}\nolimits}
\def\Ln{\mathop{\mathrm {Ln}}\nolimits}

\def\Aff{\mathop{\mathrm {Aff}}\nolimits}
\def\aff{\mathop{\mathrm {aff}}\nolimits}

\def\Id{\mathop{\mathrm {Id}}\nolimits}
\def\Ad{\mathop{\mathrm {Ad}}\nolimits}
\def\ad{\mathop{\mathrm {ad}}\nolimits}

\def\Id{\mathop{\mathrm {Id}}\nolimits}
\def\Aut{\mathop{\mathrm {Aut}}\nolimits}

\def\sgn{\mathop{\mathrm {sgn}}\nolimits}

\newtheorem{theorem}{Theorem}[section]
\newtheorem{lemma}[theorem]{Lemma}
\newtheorem{proposition}[theorem]{Proposition}
\newtheorem{corollary}[theorem]{Corollary}
\theoremstyle{remark}
\newtheorem{remark}[theorem]{Remark}

\theoremstyle{definition}
\newtheorem{definition}[theorem]{Definition}

\begin{document}
\title{Quantum Strata of Coadjoint Orbits}
\author{Do Ngoc Diep}
\address{Institute of Mathematics\\  National Centre for Science and Technology\\  P. O. Box 631, Bo Ho, VN-10000, Hanoi, Vietnam}
\email{dndiep@hn.vnn.vn} 
\address{{\rm and}\\ Department of Mathematics, The University of Iowa, 14 MacLean Hall, Iowa City, IA 52242-1419, U.S.A.}
\email{ndiep@math.uiowa.edu}

\thanks{Version of May 18, 2000.\\
This work is supported in part by the  Alexander von Humboldt Foundation, Germany and the National Foundation for Research in Fundamental Sciences, Vietnam.}
\maketitle
\begin{abstract}
In this paper we construct quantum analogs of strata of coadjoint orbits and describe their representations. This kind objects play an important role in describing quantum groups as  repeated extensions of quantum strata.
\end{abstract}

\section*{Introduction}
For locally compact groups, their C*-algebras contain exhausted informations about the groups them-selves and their representations, see \cite{diep1}, \cite{diep2}. In some sense \cite{rieffel1}-\cite{rieffel2}, the group algebras can be considered as C*-algebraic deformation quantization $C^*_q(G)$ at the special value $q=1$.

In \cite{diep1} and \cite{diep2}, the group C*-algebras were described as repeated extensions of C*-algebras of strata of coadjoint orbits. Quantum groups are {\it group Hopf algebras}, i.e. replace C*-algebras by special Hopf algebras ``of functions''. It is therefore interesting to ask whether we could describe quantum groups as some repeated extensions of some kind quantum strata of coadjoint orbits? We are attempting to give a positive answer to this question. It is not yet completely described but we obtained a reasonable answer. Let us describe the main ingredients of our approach.

For the {\it good strata}, namely families of with some good enough parameter space, of coadjoint orbits, there exist always continuous fields of polarizations (in the sense of the representation theory), satisfying {\it the L. Pukanszky irreducibility condition}: for each orbit ${\mathcal O}$ and any point $F_{\mathcal O}$ in it, the affine subspace, orthogonal to some polarizations with respect to the symplectic form  is included in orbits themselves,  i.e.
$$F_{\mathcal O} + \mathfrak h_\mathcal{O}^\perp \subseteq \mathcal{O}$$
and $$\dim \mathfrak{h}_{\mathcal O} = \frac{1}{2}\dim \mathcal{O}.$$ We choose the the canonical Darboux coordinates with impulse $p$'s-coordinates, following a vector structure basis of $\mathfrak h^\perp$. From this we can deduce that in this kind of Darboux coordinates, the Kirillov form $\omega_{\mathcal O}$ locally are canonical and every element $X\in \mathfrak g = \Lie G$ can be considered as a function $\tilde{X}$ on $\mathcal O$, linear on $p$'s-coordinates, i.e. $$\tilde{X} = \sum_{i=1}^n a_i(q)p_i + a_0(q).$$ 

This essential fact gives us a possibility to effectively write out the corresponding $\star$-product of functions,  define quantum strata $C_q(V,\mathcal F)$. On the strata acts our Lie group of symmetry. It induces therefore an action on equivariant differential operators. Using the indicate fields of polarizations we prove some kind of Poincar\'e-Birkhoff-Witt theorem and then provide quantization with separation of variable in sense of Karabegov \cite{karabegov1}. We can then express the corresponding representations of the quantum strata $C_q(V,\mathcal F)$, where $V \subseteq \bigcup_{\dim \mathcal O = const} \mathcal O$, through the Feynman path integrals etc...,see \cite{diep3}. 

Our main result is Theorem \ref{mainthm} in which we construct the $\star$-product exactly, Theorem \ref{separationthm}, where the representations are obtained from Poincar\'e-Birkhoff-Witt separation of variable and $\star$-products,  and Theorem \ref{fourierthm} where we express the product through Fourier integral operators. 

Let us now describe in brief the structure of sections. In section 1, we construct a canonical local coordinate system, where the generating functions $\tilde{X}$, for $X\in \mathfrak g$ are linear in the impulse coordinates $p$'s, Theorem 1.4. In this kind of coordinates, we can in \S2 construct construct the local $\star$-product, Theorem 2.1, then globalize it into a $\Gamma$-invariant $\star$-product on the universal coverings of coadjoint orbits and then push down to the coadjoint orbits them-selves, Theorem 2.4, where $\Gamma = \pi_1(\mathcal O)$ - the fundamental group of the coadjoint orbit $\mathcal O$. IN section 3 we let this left $\star$-multiplication of functions, acting on the quantum bundle sections. It manipulates an action of the generating functions as right $G$-invariant pseudo-differential operators of the first order. We then use the universal property of the universal enveloping algebra $U(\mathfrak g)$ to construct a quantizing homomorphism $U(\mathfrak g) \to PDO_G(\mathcal O)$, Theorem 3.1. We prove then a version of the Poincar\'e-Birkhoff-Witt Theorem, associated with a complex polarization. Theorem 3.3. In section 4, we show that this kind representations are the same as those obtained from the procedure of multidimensional quantization, \cite{diep2}. In section 5, we restrict to special case of strata of coadjoint orbits appeared from some versions of solvable cases of the Gelfand-Kirillov conjecture. In those cases we can express the associated unitary representations of the groups $G$ as some oscilatting Fourier integrals, Theorem 5.5. In Section 6, all the main ideas are demonstrated in 3 series of examples.

\section{Canonical coordinates on a stratum}

Let us consider a connected and simply connected Lie group $G$ with Lie algebra $\mathfrak g$. Denote the dual to $\mathfrak g$ vector space by $\mathfrak g^*$. It is well-known that the action of $G$ on itself by conjugation $$A(g) : G \to G,$$ defined by $A(g)(h) := ghg^{-1}$ keeps the identity element $h=e$ unmoved. This induces the tangent map $\Ad(g) := A(g)_*: {\mathfrak g} = T_eG\to \mathfrak g,$ defined by $$\Ad(g)X := \frac{d}{dt}|_{t=0} A(g)\exp(tX)$$
and the co-adjoint action 
$K(g) := Ad(g^{-1})^* $ 
maps the dual space 
${\mathfrak g}^*$ 
into itself. The orbit space 
${\mathcal O}(G):= {\mathfrak g}^*/G$ 
is in general a bad topological space, namely non-Hausdorff, in general. Consider 
an arbitrary orbit 
$\Omega\in{\mathcal O}(G)$
and an element 
$F \in {\mathfrak g}^*$ in it. The stabilizer is denote by 
$G_F$, 
its connected component by $(G_F)_0$ and its Lie algebra by 
${\mathfrak g}_F := \Lie(G_F)$. It is well-known that
$$
\begin{array}{ccc}
G_F & \hookrightarrow & G\\
    &                 & \Big\downarrow\\
    &                  & \Omega_F
    \end{array}
$$
is a principal bundle with the structural group $G_F$. 
Let us fix some {\it connection in this principal bundle, } 
\index{connection on 
principal bundle} i.e. some {\it trivialization } \index{trivialization}
of this bundle.  
We want to construct representations in some cohomology spaces 
with coefficients in the sheaves of sections of some vector bundle 
associated with this principal bundles. 
It is well know  that every vector bundle 
is an induced one with respect to some representation of the structural group
in the typical fiber.
It is natural to fix some unitary representation 
$\tilde{\sigma}$ 
of $G_F$ such that its kernel contains $(G_F)_0$, the character
$\chi_F$ of the connected component of stabilizer
$$\chi_F(\exp{X}) := \exp{(2\pi\sqrt{-1}\langle F,X\rangle )}$$
and therefore the differential 
$D(\tilde{\sigma}\chi_F)=\tilde{\rho}$ 
is some representation  of the Lie algebra ${\mathfrak g}_F$. 
We suppose that the representation $D(\tilde{\rho}\chi_F)$ was extended to 
the complexification $({\mathfrak g}_F)_{\Bbb C}$. 
The whole space of all sections seems to be so large for the construction 
of irreducible unitary representations. 
One considers the invariant subspaces with the help of some so 
called {\it polarizations}, see \cite{diep1}, \cite{diep2}.

\begin{definition} We say that a triple $(\mathfrak p_\mathcal O, \rho_\mathcal O, \sigma_{0,\mathcal O})$ is a 
$(\tilde{\sigma},F)$-{\it polarization} of $\mathcal O$ iff :
\begin{enumerate}
\item[a.] $\mathfrak p_\mathcal O$ is some complex sub-algebra of the complexified $\mathfrak g_\mathbf C$, containing $\mathfrak g_{F_\mathcal O}$.
\item[b.] The sub-algebra $\mathfrak p_\mathcal O$ is invariant under the action of all the operators of type $Ad_{\mathfrak g_\mathbf C} x, x\in G_{F_\mathcal O}.$
\item[c.] The vector space $\mathfrak p_\mathcal O + \overline{\mathfrak p_\mathcal O}$ is the complexification of some real Lie sub-algebra $\mathfrak m_\mathcal O := (\mathfrak p_\mathcal O + \overline{\mathfrak p}_\mathcal O) \cap \mathfrak g.$
\item[d.] All the subgroups $M_{0, \mathcal O}$, $H_{0, \mathcal O}$, $M_\mathcal O$, $H_\mathcal O$ are closed, where by definition, $M_{0,\mathcal O}$ (resp., $H_{0, \mathcal O}$) is the connected subgroup of $G$ with the Lie algebra $\mathfrak m_\mathcal O$ (resp., $\mathfrak h_\mathcal O := \mathfrak p_\mathcal O \cap \mathfrak g$) and the semi-direct products $M:= M_{0,\mathcal O} \ltimes G_{F_\mathcal O}$, $H_\mathcal O := H_{0,\mathcal O} \ltimes G_{F_\mathcal O}$. 
\item[e.] $\sigma_{0,\mathcal O}$ is an irreducible representation of $H_{0,\mathcal O}$ in some Hilbert space $V_\mathcal O$ such that: 
(1.) the restriction $\sigma_{0,\mathcal O}\vert_{G_{F_\mathcal O} \cap H_{0,\mathcal O}}$ is some multiple of the restriction $\chi_{F_\mathcal O}.\tilde{\sigma}_\mathcal O\vert_{G_{F_\mathcal O} \cap H_{0,\mathcal O}}$, where the character $\chi_\mathcal O$ is by definition, $\chi_\mathcal O(\exp X) = \exp(2\pi\sqrt{-1}\langle F_\mathcal O, X\rangle )$;
(2.) under the action of $G_{F_\mathcal O}$ on the dual $\hat{H}_{0,\mathcal O}$, the point $\sigma_{0,\mathcal O}$ is fixed.
\item[f.] $\rho_\mathcal O$ is a representation of the complex Lie algebra $\mathfrak p_\mathcal O$ in the same $V_\mathcal O$, which satisfies the Nelson conditions for $H_{0,\mathcal O}$ and $\rho_\mathcal O \vert_{\mathfrak h_\mathcal O} \cong D\sigma_{0,\mathcal O}$.
\end{enumerate}
\end{definition}
There is a natural order in the set of all $(\tilde{\sigma}_\mathcal O, F_\mathcal O)$-polarizations by inclusion and from now on speaking about polarizations we means always the maximal ones. It is not hard to prove that the (maximal) polarizations are also the Lagrangian distributions and in particular the co-dimension of $\mathfrak h_\mathcal O$ in $\mathfrak g$ is a half of the dimension of the coadjoint orbit $\mathcal O$,  $$codim_\mathfrak g \mathfrak h_\mathcal O = \frac{1}{2}\dim \mathcal O,$$ see e.g. \cite{diep2} or \cite{kirillov}.

Let us now recall the Pukanszky condition.
\begin{definition}
We say that the $(\tilde{\sigma}_\mathcal O,F_\mathcal O)$-polarization $(\mathfrak p_\mathcal O, \rho_\mathcal O, \sigma_{0,\mathcal O})$ satisfies the Pukanszky condition, iff $$F_\mathcal O + \mathfrak h_\mathcal O^\perp \subset \mathcal O.$$
\end{definition}

\begin{remark}
\begin{itemize}
\item Pukanszky conditions involve an inclusion of the Lagrangian affine subspace of $p$'s coordinates into the local Darboux coordinates.
\item The partial complex structure on orbits let us to use smaller subspaces of section in the induction construction, as subspaces of partially invariant, partially holomorphic sections of the induced bundles.
\end{itemize}
\end{remark}

\begin{theorem}
There exists on each coadjoint orbit a local canonical system of Darboux coordinates, in which the Hamiltonian function $\tilde{X}$, $X \in \mathfrak g$, are linear on $p$'s impulsion coordinates and in theses coordinates, $$\tilde{X} = \sum_{i=1}^n a_i(q)p_i + a_0(q).$$
\end{theorem}
\begin{proof}
We supposed that our coadjoint orbit admit at least a polarization $(\mathfrak p,\rho, \sigma)$, satisfying L. Pukanszky's condition of irreducibility
$$F_\mathcal O + \mathfrak h_\mathcal O^\perp \subseteq \mathcal O .$$ The codimension of $\mathfrak h_{\mathcal O}$ and therefore the dimension of $\mathfrak h_{\mathcal O}^\perp$ is $\frac{1}{2}\dim {\mathcal O}$. 
Consider a canonical system of geodesics. The geodesics corresponding to the affine subspace $\mathfrak h_{\mathcal O}^\perp$ provide linear coordinates $p_1, p_2, ...,p_n$. The others are the corresponding $q^1, \dots , q^n$.
Therefore we can arrange a local system of coordinates, such that exponential map gives linear geodesics on $p$'s directions. 
\end{proof}

\section{Poisson structure on strata}

Let us first recall the construction of strata of coadjoint orbit from
\cite{diep1}-\cite{diep2}. The orbit space $\mathcal O(G)$ is a disjoint union of $\Omega_{2n}$, each of which is a union of the coadjoint orbits of dimension $2n$, $ 0 \leq 2n \leq \dim G$. Denote $$V_{2n} = \cup_{\dim \mathcal O = 2n} \mathcal O.$$ Then $V_{2n}$ is the set of points of fixed rank $2n$ of the Poisson structure bilinear function $$\{X,Y\}(F) = \langle F,[X,Y]\rangle .$$

Suppose that it is a foliation, at least it is true for $V_{2n}$, where $2n$ is the maximal dimension possible in $\mathcal O(G)$. It can be shown that the foliation $V_{2n}$ can be obtained from the group action of $\mathbf R^{2n}$ on $V_{2n}$. For this aim, let us consider a basis $X_1, \dots, X_{2n}$ of the tangent space $T_{F_\mathcal O}\mathcal O \cong \mathfrak g /\mathfrak g_{F_\mathcal O}$ at the point $F_\mathcal O \in \mathcal O \subset V_{2n}$. We can define an action of $\mathbb R^{2n}$ on $V_{2n}$ as 
$$(t_1,\dots,t_{2n}) \mapsto \exp(t_1X_1)\exp(t_2X_2)\dots \exp(t_{2n}X_{2n})F_\mathcal O .$$
We have therefore the Hamiltonian vector fields
$$\xi_k := \frac{d}{dt}\vert_{t=0} \exp(t_kX_k)F, \forall k=1,\dots, 2n$$ and their span $\mathcal F_{2n} = \{\xi_1,\dots,\xi_{2n}\}$ provides a tangent distribution. It is easy to show that we have therefore a measurable (in sense of A. Connes) foliation. One can therefore define also the Connes C*-algebra $C^*(V_{2n},\mathcal F_{2n})$, $0 \leq 2n \leq \dim G$. By introducing some technical condition in \cite{diep1}, it is easy to reduce these C*-algebras to extension of other ones, those are in form of tensor product $C(X) \otimes \mathcal K(H)$ of algebras of continuous functions on compacts  and the elementary algebra $\mathcal K(H)$ of compact operators in a separable Hilbert space $H$. The strata of these kinds we means {\it good strata}. Another kind of good strata of coadjoint orbits are obtained from relation with the cases where the Gelfand-Kirillov conjecture was solved, for examples for connected and simply connected solvable Lie groups, see \S5.

It is deduced from a result of Kontsevich that this Poisson structure can be quantized. This quantization however is formal. The question of convergence of the quantizing series is not clear. We show in this section that in the case of charts with the linear p's impulse coordinates, the corresponding $\star$-product  is convergent.

In this kind of special local chart systems of Darboux coordinates it is easy to deduce existence local convergent $\star$-products. 
\begin{theorem}
Locally on each coadjoint orbit, there exist a convergent $star$-product.
\end{theorem}
\begin{proof} 
Let us denote by $\mathcal F^{-1}_p$ the Fourier inverse transformation on variables $p$'s and by $\mathcal F_p$ the Fourier transformation. Let us denote by $PDO_G(\mathcal O)$ the algebra of $G$-invariant pseudodifferential operators on $\mathcal O$. Locally, the Fourier transformation maps symbols (as function on local coordinates of $\mathcal O$) to $G$-invariant pseudodifferential operators and conversely, the inverse Fourier transformation  maps the pseudodifferential operators to some specific classes of symbols. For two symbols $f,g \in \mathcal F^{-1}_p(PDO_G(\mathcal O))= \mathbf C_q(\mathcal O)$, their Fourier images $\mathcal F_p(f)$ and $\mathcal F_p(g)$ are operators and we can define their operator product and then take the Fourier inverse transforms, as $\star$-product
$$f \star g := \mathcal F^{-1}_p(\mathcal F_p(f).\mathcal F_p(g)).$$ So this product is again a symbol n the same class $\mathbf C_q(\mathcal O)$.

Because of existence of special coordinate systems, linear on $p$'s coordinates we can treat for the good strata in the same way as in the cases of exponential groups.  And the formal series of $\star$-product is convergent.
\end{proof}

\begin{remark}
The proof could be also done in the same scheme as in exponential or compact cases, see \cite{arnalcortet1}-\cite{arnalcortet2}.
\end{remark}

Let us denote by $\Gamma= \pi_1({\mathcal O})$ the fundamental group of the orbit. Our next step is to globally extend this kind of local $\star$-products. Our idea is as follows. We lift this $\star$-product to the universal covering of coadjoint orbits as some $\Gamma$-invariant $\star$-products, globally extend them in virtue of the monodromy theorem and then pushdown to our orbits.
We start with the following lemmas

\begin{lemma} \label{lemma1}
There is one-to-one correspondence between  $\star$-products on Poisson manifolds and $\Gamma$-invariant $\star$-products on their universal coverings.
\end{lemma}
\begin{proof}
By the lifting properties of the universal covering, we can easily lift each $\star$-product on a Poisson manifold onto its universal covering. This correspondence is one-to-one.
\end{proof}

We use this lemma to describe existence of a $\star$-product on coadjoint orbits.

\begin{lemma}\label{lemma2}
On a universal covering, each local $\star$-product can be uniquely extended to some global $\star$-product on this covering. 
\end{lemma}
\begin{proof} For local charts, there exist deformation quantization, as said above, $f \mapsto Op(f)$ by using the formulas of Fedosov quantization. Also in the intersection of two local charts of coordinates $(q,p)$ and $(\tilde{q}, \tilde{p})$, there is a symplectomorphism, namely $\varphi$ such that $$(\tilde{q},\tilde{p}) = \varphi(q,p).$$ Using the local oscilatting integrals and by compensating the local Maslov's index obstacles, one can exactly construct the unitary operator $U$ such that $$Op(f\circ \varphi) = U.Op(f) .U^{-1},$$ see for example Fedosov's book \cite{fedosov}.
Because the universal coverings are simply connected, the extensions can be therefore produced because of Monodromy Theorem.
\end{proof}

\begin{remark}
The operators ``U'' of this kind, depending on two local charts as parameters, provide some cohomological 2-class and therefore are classified by some $2^{nd}$ cohomology class with values in unitary operators. This reduces to some classification of all the possibilities of quantizations, upto the first cohomology classes, i.e. upto conjugations.
\end{remark}

\begin{theorem}\label{mainthm}
There exists a convergent $\star$-product on each orbit, which is a symplectic leaf of the Poisson structure on each stratum of coadjoint orbits.
\end{theorem}
\begin{proof}
From the description of the canonical coordinates in the previous section, we see that there exists at least a convergent $\Gamma$-invariant local $\star$-product. This local $\star$-product then extended to  a $\Gamma$-invariant global one on the universal covering, which produces a convergent $\star$-product on the coadjoint orbits, following Lemmas \ref{lemma1},\ref{lemma2}.
\end{proof}

\section{Star-product, Quantization and PBW Theorem}

We use the constructed in the previous section $\star$-product to provide an action of functions on the spaces of partially invariant partially holomorphic sections of the corresponding partially invariant holomorphically induced bundles associated with polarizations. It is possibles because the quantum induced bundles are locally trivial and the spaces of partially invariant partially holomorphic sections with section in Hilbert spaces  locally are finitely generated modules over the algebras of quantizing functions, see \cite{diep2}. 
We use then the construction of Karabegov and Fedosov to obtain a Hopf *-algebra of longitudinal pseudo-differential operators elliptic along the leaves of this measurable foliations. The main ingredient is that we use here the Poincar\'e-Birkhopf-Witt theorem to provide this quantization. 

\begin{theorem}\label{theorem3.1}
There is a natural deformation quantization with separation of variables, corresponding to the Poincar\'e-Birkhoff-Witt Theorem for polarizations.
\end{theorem}

We have from the multidimensional quantization, see e.g. \cite{diep1}-\cite{diep2}, $\mathfrak g \to PDO^1(\mathcal O)$, $X \mapsto \hat{X}$. More precisely, Let us denote by $PDO(\mathcal O)$ the algebra of right $G$-invariant pseudo-differential operators on $\mathcal O$, i.e. the continuous maps from $C^\infty(\mathcal O)$ into itself not extending support. We also denote $PDO^1(\mathcal O)$ the Lie algebras of right $G$-invariant first order pseudo-differential operators. By the procedure of multidimensional quantization, there is a homomorphism of Lie algebras
$$\mathfrak g \to PDO_G^1(\mathcal O) \subset PDO_G(\mathcal O).$$ Following the universal property of $U(\mathfrak g)$, there is a unique homomorphism of associative algebras $U(\mathfrak g) \to PDO_G(\mathcal O)$making the following diagram commutative
$$\CD
\mathfrak g @>>> PDO_G(\mathcal O)\\
@VVV    @AAA\\
U(\mathfrak g) @= U(\mathfrak g)
\endCD$$

 Let us first describe the machinery applied in the method of Karabegov's separation of variable. We use the idea about polarizations in multidimensional quantization, \cite{diep2}.

Let us recall the root decomposition $$\mathfrak g = \mathfrak n_- \oplus \mathfrak a \oplus \mathfrak n_+.$$ If $\mathfrak p$ is a complex polarization, then from definition we have
$\mathfrak m_\mathbf C = \mathfrak p \oplus \bar{\mathfrak p},$
$\mathfrak h_\mathbf C = \mathfrak p \cap \bar{\mathfrak p},$
where $\mathfrak m = \mathfrak g \cap (\mathfrak p \oplus \bar{\mathfrak p}) $ and 
$\mathfrak h = \mathfrak g \cap (\mathfrak p \oplus \bar{\mathfrak p}).$ The quotients subspaces $\mathfrak p / \mathfrak h_{\mathbf C}$ and $\bar{\mathfrak p}/\mathfrak h_{\mathbf C}$ are included in $\mathfrak m$ as linear subspaces (not necessarily to be sub-algebras). Let us fix some (non-canonical) inclusions.
Let us denote by $U(\mathfrak p/\mathfrak h_{\mathbf C})$ (resp. $U(\bar{\mathfrak p}/\mathfrak h_{\mathbf C})$) the {\it sub-algebra, generated} by elements from $\mathfrak p/\mathfrak h_{\mathbf C} \hookrightarrow \mathfrak m_{\mathbf C}$ (resp. $\mathfrak p/\mathfrak h_{\mathbf C} \hookrightarrow \mathfrak m_{\mathbf C}$) in the universal algebra $U(\mathfrak m)$. We have therefore an analog of the well-known Poincar\'e-Birkhoff-Witt theorem.
\begin{theorem}[Poincar\'e-Birkhoff-Witt Theorem]
If  $\mathfrak p$ is as polarization, then $$U(\mathfrak m_{\mathbf C}) \cong U(\mathfrak p/\mathfrak h_{\mathbf C}) \otimes U(\mathfrak h_{\mathbf C}) \otimes U(\bar{\mathfrak p} /\mathfrak h_{\mathbf C})$$ 
\end{theorem}
\begin{proof}
Let us fix in a basis each of $\mathfrak p/\mathfrak h_\mathbf C$, $\mathfrak h_\mathbf C$ and $\bar{\mathfrak p}/\mathfrak h_\mathbf C$. We have therefore a basis of $\mathfrak m_\mathbf C = \mathfrak p /\mathfrak h_\mathbf C \oplus \mathfrak h_\mathbf C \oplus \bar{\mathfrak p}/\mathfrak h_\mathbf C$. Our theorem is therefore deduced from the original Poincar\'e-Birkhoff-Witt Theorem \cite{pbw}.
\end{proof}

It is easy to see that by this reason, on the variety $M/H$ there is a natural complex structure and therefore on coadjoint orbits exists some partial complex structure. We use this complex structure and this Poincar\'e-Birkhoff-Witt theorem to do a separation of variables on $M$ and apply the machinery of Karabegov.

\begin{remark}
Because of BKW, $U(\mathfrak g) \cong U(\mathfrak p_{\mathcal O}/\mathfrak h_{\mathcal O}) \otimes U(\mathfrak h_{\mathcal O}) \otimes U(\overline{\mathfrak p}_{\mathcal O}/\mathfrak h_{\mathcal O})$ and because of this our quantizing map is coincided with that one used by Karabegov in the quantization with separation of variables in case of coadjoint orbits with totally complex polarizations $\mathfrak g_\mathbf C = \mathfrak m_\mathbf C = \mathfrak p \oplus \overline{\mathfrak p}$.
\end{remark}

\begin{theorem} \label{separationthm}
In the case of totally complex polarizable coadjoint orbits, the quantization map from the theorem \ref{theorem3.1} is coincided with the rule of Karabegov's quantization with separation of variables.
\end{theorem}

\section{Representations}

Let us consider now a continuous fields of $(\tilde{\sigma}_\mathcal O,F_\mathcal O)$-polarizations $(\mathfrak p_\mathcal O, \rho_\mathcal O, \sigma_{0,\mathcal O})$ satisfying the Pukanszky condition. On one hand side, we can use the multidimensional quantization procedure to obtain irreducible unitary representations $\Pi_\mathcal O =  Ind(G;\mathfrak p_\mathcal O,H_\mathcal O,\rho_\mathcal O,\sigma_{0,\mathcal O})$ of $G$, \cite{diep2}. On other hand side, we can use $\star$-product construction to provide the representations $\Pi_\mathcal O : \tilde{X\vert_\Omega} \mapsto \hat{\ell}_X$ of the quantum strata $\mathbf C_q(\Omega)$. We'd like to show we have the same one.

\begin{definition}
{\it Quantum coadjoint orbit} $\mathbf C_q(\mathcal O)$ is defined as the Hopf algebra of symbols of differential operators $U_{q,\mathcal O}(\mathfrak g) = U(\mathfrak g)\vert_{\mathcal O}$. The homomorphism $Q: U(\mathfrak g) \to PDO_G(\mathcal O)$ is defined to be the {\it second quantization homomorphism}.
\end{definition}

\begin{theorem}
 The representation of the Lie algebra obtained from $\star$-product is equal to the representations obtained from the multidimensional quantization procedure.
\end{theorem}
\begin{proof} Let us recall \cite{diep2}, that 
$$\Lie_X \; Ind(G;\mathfrak p_\mathcal O,H_\mathcal O,\rho_\mathcal O,\sigma_{0,\mathcal O}) \cong \hat{X}.$$
From the construction of quantization map $$U(\mathfrak g) \to PDO_G(\mathcal O)$$ as the map arising from the universal property the map $\mathfrak g \to PDO^1_G(\mathcal O),$ $$\hat{\ell}_X = Op(\tilde{X}) = \hat{X}, \forall X \in \mathfrak g = \Lie(G).$$ The associated representation of $\mathbf C_q(\Omega)$ is the solution of the Cauchy problem for the differential equation
$$\frac{\partial}{\partial t}U(t, q,p) = \ell_X U(t,q,p),$$
$$U(0,q,p) = \Id .$$
The solution of this problem is uniquely defined.
\end{proof}

\section{Oscillating Fourier integrals}

Let us in this section consider the good family of coadjoint orbits arising from the solved cases of Gelfand-Kirillov Conjecture. The results in this section is revised from the \cite{diep3}.

Consider a connected and simply connected Lie group $G$ with Lie algebra $\mathfrak g$.
\begin{theorem}
There exists a $G$-invariant Zariski open set $\Omega$ and a covering $\tilde{\Omega}$ of $\Omega$, with natural action of $G$ such that for each continuous field of polarizations $(\mathfrak p_\mathcal O, H_\mathcal O, \rho_\mathcal O, \sigma_{0,\mathcal O})$, $\mathcal O \in \Omega/G$, the Lie derivative of the direct integral of representations arized from the multidimensional quantization procedure is equivalent to the tensor product of the Schr\"odinger representation $\Pi = Sch$ of the Gelfand-Kirillov basis of the enveloping field and a continuous field of trivial representations $\{V_\mathcal O\}$ on $\tilde{\Omega}$.
\end{theorem}
\begin{proof}
Our proof is rather long and consists of several steps:

1. We apply the construction of Nghiem \cite{nghiem5} to the solvable radical $^r{\mathfrak g}$ of $\mathfrak g$ to obtain the Zariski open set $^r\Omega$ in $^r\mathfrak g^*$ and its covering $^r\tilde{\Omega}$.

2. The general case is reduced to the semi-simple case $^s\mathfrak g$, see (\cite{nghien1}, Thm. C). Denote by $^sG$ the corresponding analytic subgroup of $G$.

3. Take the Zariski open set $\mathcal A_s\mathcal P(^sG)$ of admissible and well-polarizable strongly regular functionals from $^s\mathfrak g^*$ and its covering $\mathcal B_s(^sG)$ via Duflo's construction \cite{duflo1}.

4. The desired $G$-invariant Zariski open set and its covering are the corresponding Cartesian products
$$\Omega= \mathcal A_s\mathcal P(^sG) \times {}^r\Omega^0 \times {}^r\Omega^1 \times \dots \times {}^r\Omega^k$$
$$\tilde{\Omega}= \mathcal B_s(^sG) \times {}^r\tilde{\Omega^0} \times {}^r\tilde{\Omega^1} \times \dots\times  {}^r\tilde{\Omega}^k$$

5. \begin{lemma}
There exists a continuous field of polarizations of type $(\mathfrak p_\mathcal O, H_\mathcal O,\rho_\mathcal O,\sigma_{0,\mathcal O})$, for each $\mathcal O \in \tilde{\Omega}/G$.
\end{lemma}

6. \begin{lemma}
The Lie derivative commutes with direct integrals, i.e. $$\nabla \int^\oplus_{\tilde{\Omega}/G}Ind(G;\mathfrak p_\mathcal O, H_\mathcal O,\rho_\mathcal O,\sigma_{0,\mathcal O})d\mathcal O = \int^\oplus_{\tilde{\Omega}/G} \nabla\;Ind(G;\mathfrak p_\mathcal O, H_\mathcal O,\rho_\mathcal O, \sigma_{0,\mathcal O}) d\mathcal O.$$
\end{lemma}

7. \begin{lemma}
The restriction of the Schr\"odinger representation to coadjoint orbits provides a continuous field of polarizations. In particular,
$$D\sigma_{0,\mathcal O} = \rho_\mathcal O\vert_{\mathfrak p_{\mathcal O} \cap \mathfrak g} \cong mult \; Sch_\mathcal O.$$
\end{lemma}

\end{proof}

\begin{theorem}\label{fourierthm}
1. There exists an operator-valued phase $\Phi(t,z,x)$ and an operator-valued amplitude $a(t,z,x,y,\xi)$ extended from the expression
$$\exp(\Phi(t,z,y) + \sqrt{-1}\xi(g(t)x-x))$$ in such a fashion that the action of the representation $$\Pi_\mathcal O = Ind(G;\mathfrak p_\mathcal O, H_\mathcal O,\rho_\mathcal O,\sigma_{0,\mathcal O})$$ can be expressed as an oscilatting Fourier integral
$$\Pi_\mathcal O(g(t))f(z,x) = c.\int_{\mathbf R^M}\int_{\mathbf R^M} a(t,z,x,y,\xi)\exp(\sqrt{-1}\xi(x-y))f(z,y)dyd\xi,$$ where $c$ is some constant.

2. For each function $\varphi$ of Schwartz class $\mathcal S(G)$ satisfying the compactness criteria \cite{diep2} in every induction step \cite{lipsman}, \cite{diep2}, the operator $\Pi_\mathcal O(\varphi)$ is of trace class and its action can be expressed as the oscilatting Fourier integral 
$$\Pi_\mathcal O(\varphi)f(z,x) = const.\int_{\mathbf R^M}\int_{\mathbf R^M}(\int_{\mathbf R^\ell} a(t,z,x,y,\xi)\varphi(g(t))dt) \times$$ $$\times \exp(\sqrt{-1}\xi(x-y))f(z,y)dyd\xi. $$
Hence, its trace is $$tr\Pi_\mathcal O(\varphi) = const.\int_{\mathbf R^M}\int_{\mathbf R^M}(\int_{\mathbf R^\ell} a(t,z,x,x,\xi)\varphi(g(t))dt )dxd\xi$$
\end{theorem}
\begin{proof}

The proof also consists of several steps:

1. From \cite{diep1} - \cite{diep2} and \cite{duflo1} it is easy to see obtain a slight unipotent modification( i.e. a reduction to the unipotent radical) of the multidimensional quantization procedure. We refer the reader to \cite{tranvui1}-\cite{tranvui2}, and \cite{trandong} for a detailed exposition.

2. From the unitary representations $\Pi_\mathcal O = Ind(G; \mathfrak p_\mathcal O, H_\mathcal O, \rho_\mathcal O , \sigma_{0,\mathcal O})$ in the unipotent context and its differential $\pi_\mathcal O = D\Pi_\mathcal O(G;  \mathfrak p_\mathcal O, H_\mathcal O, \rho_\mathcal O, \sigma_{0,\mathcal O})$, it is easy to select the constant term and the vector fields term, (see \cite{nghiem4} for the simplest case and notations), 
$$\pi(L) = \pi_c(L) + \pi_d(L), \quad \pi_s(L) = a_0(L, x),$$
$$\pi_d(L) = \partial (L) = \sum_k a_k(L) \frac{\partial}{\partial x_k} .$$

The operator-valued partial (i.e. depending on $L$) phase $$\Phi_L(t_L,z,x) = \int_0^{t_L} a_0(z,g(s)x)ds.$$ The phase $\Phi(t,z,x)$ is the sum of the partial phases, on which the one-parameter group $g_L(s)$ operates by translations for all the factors on the left of $g_L(t_L)$ in the ordered product $$g(t) = \prod_L g_L(t_L). $$

It is easy then to see that our induced representations $\Pi_\mathcal O$ act as follows
$$\Pi_\mathcal O(g(t))f(z,x) = \sigma_{0,\mathcal O}(t,z,x) f(z,g(t)x).$$

3. It suffices now to apply the Fourier transforms
$$f(z,x) = (2\pi )^M \int_{\mathbf R^M} \int_{\mathbf R^M} \exp\{\sqrt{-1}\xi (x-y)\} f(z,y)dyd\xi$$ to the function $\Pi_\mathcal O(g(t))f(z,x)$ to obtain
$$\Pi_\mathcal O(g(t))f(z,x) = $$
$$(2\pi)^M\int_{\mathbf R^M}\int_{\mathbf R^M} a(t,z,x,y,\xi)\exp\{\sqrt{-1}\xi (x-y)\}f(z,y)dyd\xi,$$
where the amplitude $a(t,z,x,y,\xi)$ is the natural extension of the expression $$\exp\{ \Phi(t,z,y) + \sqrt{-1}\xi (g(t)x -x) \}$$ in the correspondence with the fields of polarizations. Hence for each $\varphi \in \mathcal S(G)$,
$$\Pi_\mathcal O(\varphi)f(z,x) = (2\pi)^M \int_{\mathbf R^M} \int_{\mathbf R^M} (\int_{\mathbf R^\ell} a(t,z,x,y\xi)\varphi(g(t))dt) \times$$
$$\times \exp\{\sqrt{-1}\xi (x - y) \}f(z,y)dyd\xi .$$

Remark that the integral $\int_{\mathbf R^\ell} a(t,z, x,y,\xi)\varphi(g(t)) dt$ is just a type of Feynman path integrals.

4. \begin{lemma}
If in every repeated induction step, see \cite{lipsman}, \cite{diep2}, $\Pi_\mathcal O(\varphi)$ satisfies the compactness criteria, then the  operator $\Pi_\mathcal O(\varphi)$ is trace class and hence $$tr\; \Pi_\mathcal O(\varphi) = (2\pi)^M \int_{\mathbf R^M}\int_{\mathbf R^M} tr( \int_{\mathbf R^\ell} a(t,z,x,x,\xi)\varphi(g(t))dtdxd\xi.$$
\end{lemma}
The proof of the theorem is therefore achieved.
\end{proof}

\section{Examples}
In this section, we illustrate the main ideas in examples.
\subsection{$\overline{MD}$-groups}
We expose in this subsection our joint works with Nguyen Viet Hai \cite{diephai}-\cite{diephai2}.
\subsubsection{The group of affine transformations of the real straight line}
We refer the reader to the work \cite{diephai} for a detailed exposition with complete proof.

{\it Canonical coordinates on the upper half-planes.}
Recall that the Lie algebra $\mathfrak g = \aff(\mathbf  R)$ of affine transformations of the real straight line is described as follows, see for example \cite{diep2}: The Lie group $\Aff(\mathbf R)$ of affine transformations of type $$x \in \mathbf R \mapsto ax + b, \mbox{ for some parameters }a, b \in \mathbf R, a \ne 0.$$ It is well-known that this group $\Aff(\mathbf R)$ is a two dimensional Lie group which is isomorphic to the group of matrices
$$\Aff(\mathbf R) \cong \{\left (\begin{array}{cc} a & b \\ 0 & 1 \end{array} \right) \vert a,b \in \mathbf R , a \ne 0 \}.$$ We consider its connected component $$G= \Aff_0(\mathbf R)= \{\left (\begin{array}{cc} a & b \\ 0 & 1 \end{array} \right) \vert a,b \in \mathbf R, a > 0 \}$$ of identity element. Its Lie algebra is
$$\mathfrak g = \aff(\mathbf R) \cong  \{\left (\begin{array}{cc} \alpha & \beta \\ 0 & 0 \end{array} \right) \vert \alpha, \beta  \in \mathbf R \}$$  admits a basis of two generators $X, Y$ with the only nonzero Lie bracket $[X,Y] = Y$, i.e. 
$$\mathfrak g = \aff(\mathbf R) \cong \{ \alpha X + \beta Y \vert [X,Y] = Y, \alpha, \beta \in \mathbf R \}.$$
The co-adjoint action of $G$ on $\mathfrak g^*$ is given (see e.g. \cite{arnalcortet2}, \cite{kirillov1}) by $$\langle K(g)F, Z \rangle = \langle F, \Ad(g^{-1})Z \rangle, \forall F \in \mathfrak g^*, g \in G \mbox{ and } Z \in \mathfrak g.$$ Denote the co-adjoint orbit of $G$ in $\mathfrak g$, passing through $F$ by 
$$\Omega_F = K(G)F :=  \{K(g)F \vert F \in G \}.$$ Because the group $G = \Aff_0(\mathbf R)$ is exponential (see \cite{diep2}), for $F \in \mathfrak g^* = \aff(\mathbf R)^*$, we have 
$$\Omega_F = \{ K(\exp(U)F | U \in \aff(\mathbf R) \}.$$
It is easy to see that
$$\langle K(\exp U)F, Z \rangle = \langle F, \exp(-\ad_U)Z \rangle.$$ It is easy therefore to see that
$$K(\exp U)F = \langle F, \exp(-\ad_U)X\rangle X^*+\langle F, \exp(-\ad_U)Y\rangle Y^*.$$
For a general element $U = \alpha X + \beta Y \in \mathfrak g$, we have
$$\exp(-\ad_U) = \sum_{n=0}^\infty \frac{1}{n!}\left(\begin{array}{cc}0 & 0 \\ \beta & -\alpha \end{array}\right)^n = \left( \begin{array}{cc} 1 & 0 \\ L & e^{-\alpha} \end{array} \right),$$ where $L = \alpha + \beta + \frac{\alpha}{\beta}(1-e^\beta)$. This means that
$$K(\exp U)F = (\lambda + \mu L) X^* + (\mu e^{\-\alpha})Y^*. $$ From this formula one deduces  \cite{diep2} the following description of all co-adjoint orbits of $G$ in $\mathfrak g^*$:
\begin{itemize}
\item If $\mu = 0$, each point $(x=\lambda , y =0)$ on the abscissa ordinate corresponds to a 0-dimensional co-adjoint orbit $$\Omega_\lambda = \{\lambda X^* \}, \quad \lambda \in \mathbf R .$$
\item For $\mu \ne 0$, there are two 2-dimensional co-adjoint orbits: the upper half-plane $\{(\lambda , \mu) \quad\vert\quad \lambda ,\mu\in \mathbf R , \mu > 0 \}$ corresponds to the co-adjoint orbit
\begin{equation} \Omega_{+} := \{ F = (\lambda + \mu L)X^* + (\mu e^{-\alpha})Y^* \quad \vert \quad \mu > 0 \}, \end{equation}
and the lower half-plane $\{(\lambda , \mu) \quad\vert\quad \lambda ,\mu\in \mathbf R , \mu < 0\}$ corresponds to the co-adjoint orbit
\begin{equation} \Omega_{-} := \{ F = (\lambda + \mu L)X^* + (\mu e^{-\alpha})Y^* \quad \vert \quad \mu < 0 \}. \end{equation}
\end{itemize}
We shall work from now on for the fixed co-adjoint orbit $\Omega_+$. The case of the co-adjoint orbit $\Omega_-$ is similarly treated. First we study the geometry of this orbit and introduce some canonical coordinates in it.
It is well-known from the orbit method \cite{kirillov} that the Lie algebra $\mathfrak g = \aff(\mathbf R)$, realized by the complete right-invariant Hamiltonian vector fields on co-adjoint orbits $\Omega_F \cong G_F \setminus G$ with flat (co-adjoint) action of the Lie group $G = \Aff_0(\mathbf R)$. On the orbit $\Omega_+$ we choose a fix point $F=Y^*$. It is well-known from the orbit method that we can choose an arbitrary point $F$ on $\Omega_F$. It is easy to see that the stabilizer of this (and therefore of any) point  is trivial $G_F = \{e\}$. We identify therefore $G$ with $G_{Y^*}\setminus G$. There is a natural diffeomorphism $\Id_{\mathbf R} \times \exp(.)$ from the standard symplectic space $\mathbf R^2$ with symplectic 2-form $dp \wedge dq$ in canonical Darboux $(p,q)$-coordinates, onto the upper half-plane $\mathbf H_+ \cong \mathbf R \rtimes
 \mathbf R_+$ with coordinates $(p, e^q)$, which is, from the above coordinate description, also diffeomorphic to the co-adjoint orbit $\Omega_+$. We can use therefore $(p,q)$ as the standard canonical Darboux coordinates in $\Omega_{Y^*}$. There are also non-canonical Darboux coordinates $(x,y) = (p,e^q)$ on $\Omega_{Y^*}$. We show now that in these coordinates $(x,y)$, the Kirillov form looks like $\omega_{Y^*}(x,y) = \frac{1}{y}dx \wedge dy$, but in the canonical Darboux coordinates $(p,q)$, the Kirillov form is just the standard symplectic form $dp \wedge dq$. This means that there are  symplectomorphisms between the standard symplectic space $\mathbf R^2, dp \wedge dq)$, the upper half-plane $(\mathbf H_+, \frac{1}{y}dx \wedge dy)$ and the co-adjoint orbit $(\Omega_{Y^*},\omega_{Y^*})$. 
Each element $Z\in \mathfrak g$ can be considered as a linear functional $\tilde{Z}$ on co-adjoint orbits, as subsets of $\mathfrak g^*$, $\tilde{Z}(F) :=\langle F,Z\rangle$.  It is well-known that this linear function is just the Hamiltonian function associated with the Hamiltonian vector field $\xi_Z$, which represents $Z\in \mathfrak g$ following the formula 
$$(\xi_Zf)(x) := \frac{d}{dt}f(x\exp (tZ))|_{t=0}, \forall f \in C^\infty(\Omega_+).$$ 
The Kirillov form $\omega_F$ is defined by the formula 
\begin{equation}\label{7} \omega_F(\xi_Z,\xi_T) = \langle F,[Z,T]\rangle, \forall Z,T \in \mathfrak g = \aff(\mathbf R). \end{equation} This form defines the symplectic structure and the Poisson brackets on the co-adjoint orbit $\Omega_+$. For the derivative along the direction $\xi_Z$ and the Poisson bracket we have relation $\xi_Z(f) = \{\tilde{Z},f\}, \forall f \in C^\infty(\Omega_+)$. It is well-known in differential geometry that the correspondence
$Z \mapsto \xi_Z, Z \in \mathfrak g$ defines a representation of our Lie algebra by vector fields on co-adjoint orbits. If the action of $G$ on $\Omega_+$ is flat \cite{diep2}, we have the second Lie algebra homomorphism from  strictly Hamiltonian right-invariant vector fields into the Lie algebra of smooth functions on the orbit with respect to the associated Poisson brackets.

Denote by $\psi$ the indicated symplectomorphism from $\mathbf R^2$ onto $\Omega_+$
$$(p,q) \in \mathbf R^2 \mapsto \psi(p,q):= (p,e^q) \in \Omega_+$$
\begin{proposition}
1. Hamiltonian function $f_Z = \tilde{Z}$ in canonical coordinates $(p,q)$ of the orbit $\Omega_+$ is of the form $$\tilde{Z}\circ\psi(p,q) = \alpha p + \beta e^q, \mbox{ if  } Z = \left(\begin{array}{cc} \alpha & \beta \\ 0 & 0 \end{array} \right).$$

2. In the canonical coordinates $(p,q)$ of the orbit $\Omega_+$, the Kirillov form $\omega_{Y^*}$ is just the standard form $\omega = dp \wedge dq$.
\end{proposition}

{\it Computation of generators $\hat{\ell}_Z$}
Let us denote by $\Lambda$ the 2-tensor associated with the Kirillov standard form $\omega = dp \wedge dq$ in canonical Darboux coordinates. We use also the multi-index notation. Let us consider the well-known Moyal $\star$-product of two smooth functions $u,v \in C^\infty(\mathbf R^2)$, defined by
$$u \star v = u.v + \sum_{r \geq 1} \frac{1}{r!}(\frac{1}{2i})^r P^r(u,v),$$ where
$$P^r(u,v) := \Lambda^{i_1j_1}\Lambda^{i_2j_2}\dots \Lambda^{i_rj_r}\partial_{i_1i_2\dots i_r} u \partial_{j_1j_2\dots j_r}v,$$ with $$\partial_{i_1i_2\dots i_r} := \frac{\partial^r}{\partial x^{i_1}\dots \partial x^{i_r}}, x:= (p,q) = (p_1,\dots,p_n,q^1,\dots,q^n)$$ as multi-index notation. It is well-known that this series converges in the Schwartz distribution spaces $\mathcal S (\mathbf R^n)$. We apply this to the special case $n=1$. In our case we have only $x = (x^1,x^2) = (p,q)$. 
\begin{proposition}\label{3.1}
In the above mentioned canonical Darboux coordinates $(p,q)$ on the orbit $\Omega_+$, the Moyl $\star$-product satisfies the relation
$$i\tilde{Z} \star i\tilde{T} - i\tilde{T} \star i\tilde{Z} = i\widetilde{[Z,T]}, \forall Z, T \in \aff(\mathbf R).$$
\end{proposition}

Consequently, to each adapted chart $\psi$ in the sense of \cite{arnalcortet2}, we associate a $G$-covariant $\star$-product.

\begin{proposition}[see \cite{gutt}]
Let $\star$ be a formal differentiable $\star$-product on $C^\infty(M, \mathbf R)$, which is covariant under $G$. Then there exists a representation $\tau$ of $G$ in $\Aut N[[\nu]]$ such that 
$$\tau(g)(u \star v) = \tau(g)u \star \tau(g)v.$$
\end{proposition}

Let us denote by $\mathcal F_pu$ the partial Fourier transform \cite{meisevogt} of the function $u$ from the variable $p$ to the variable $x$, i.e.
$$\mathcal F_p(u)(x,q) := \frac{1}{\sqrt{2\pi}}\int_{\mathbf R} e^{-ipx} u(p,q)dp.$$ Let us denote by $ \mathcal F^{-1}_p(u) (x,q)$ the inverse Fourier transform. 
\begin{lemma}\label{lem3.1}
1. $\partial_p \mathcal F^{-1}_p(p.u) = i \mathcal F^{-1}_p(x.u)  $ ,

2. $ \mathcal F_p(v) = i \partial_x\mathcal F_p(v)  $ ,

3. $P^k(\tilde{Z},\mathcal F^{-1}_p(u)) = (-1)^k \beta e^q \frac{\partial^k\mathcal F^{-1}_p(u)}{\partial^kp}, \mbox{ with } k \geq 2.$
\end{lemma}

For each $Z \in \aff(\mathbf R)$, the corresponding Hamiltonian function is $\tilde{Z} =  \alpha p + \beta e^q $ and we can consider the operator $\ell_Z$ acting on dense subspace $L^2(\mathbf R^2, \frac{dpdq}{2\pi})^\infty$ of smooth functions by left $\star$-multiplication by $i \tilde{Z}$, i.e. $\ell_Z(u) = i\tilde{Z} \star u$. It is then continuated to the whole space $L^2(\mathbf R^2, \frac{dpdq}{2\pi})$. It is easy to see that, because of the relation in Proposition (\ref{3.1}), the correspondence $Z \in \aff(\mathbf R) \mapsto \ell_Z = i\tilde{Z} \star .$ is a representation of the Lie algebra $\aff(\mathbf R)$ on the space $N[[\frac{i}{2}]]$ of formal power series in the parameter $\nu = \frac{i}{2}$ with coefficients in $N = C^\infty(M,\mathbf R)$, see e.g. \cite{gutt} for more detail.

We study now the convergence of the formal power series. In order to do this, we look at the $\star$-product of $i\tilde{Z}$ as the $\star$-product of symbols and define the differential operators corresponding to $i\tilde{Z}$. It is easy to see that the resulting correspondence is a representation of $\mathfrak g $ by pseudo-differential operators. 

\begin{proposition}
For each $Z \in \aff(\mathbf R)$ and for each compactly supported $C^\infty$ function $u \in C^\infty_0(\mathbf R^2)$, we have 
$$\hat{\ell}_Z(u) := \mathcal F_p \circ \ell_Z \circ \mathcal F^{-1}_p(u) = \alpha (\frac{1}{2}\partial_q - \partial_x)u + i\beta e^{q -\frac{x}{2}}u.$$
\end{proposition}

\begin{remark}{\rm
Setting new variables $s = q - \frac{x}{2}$, $t = q + \frac{x}{2}$, we have 
\begin{equation}
\hat{\ell}_Z(u) = \alpha\frac{\partial u}{\partial s} + i\beta e^s u,
\end{equation}
e.i. $$\hat{\ell}_Z = \alpha\frac{\partial }{\partial s} + i\beta e^s ,$$ which provides a representation of the Lie algebra $\aff (\mathbf R)$. 
}\end{remark}

{\it The associate irreducible unitary representations}

Our aim in this section is to exponentiate the obtained representation $\hat{\ell}_Z$ of the Lie algebra $\aff(\mathbf R)$ to the corresponding representation of the Lie group $\Aff_0(\mathbf R)$. We shall prove that the result is exactly the irreducible unitary representation $T_{\Omega_+}$ obtained from the orbit method or Mackey small subgroup method applied to this group $\Aff(\mathbf R)$.
Let us recall first the well-known list of all the irreducible unitary representations of the group of affine transformation of the real straight line.
\begin{theorem} [\cite{gelfandnaimark}]\label{4.1}
Every irreducible unitary representation of the group $\Aff(\mathbf R)$ of all the affine transformations of the real straight line, up to unitary equivalence, is equivalent to one of the pairwise nonequivalent representations:
\begin{itemize} 
\item the infinite dimensional representation $S$, realized in the space $L^2(\mathbf R^*, \frac{dy}{\vert y\vert})$, where $\mathbf R^* = \mathbf R \setminus \{0\}$ and is defined by the formula
$$(S(g)f)(y) := e^{iby}f(ay), \mbox{ where } g = \left(\begin{array}{cc} a & b\\ 0 & 1 \end{array}\right),$$
\item the representation $U^\varepsilon_\lambda$, where $\varepsilon = 0,1$, $\lambda \in \mathbf R$, realized in the 1-dimensional Hilbert space $\mathbf C^1$ and is given by the formula
$$U^\varepsilon_\lambda(g) = \vert a \vert^{i\lambda}(\sgn a)^\varepsilon .$$
\end{itemize}
\end{theorem}
Let us consider now the connected component $G= \Aff_0(\mathbf R)$. The irreducible unitary representations can be obtained easily from the orbit method machinery.
\begin{theorem}
The representation $\exp(\hat{\ell}_Z)$ of the group $G=\Aff_0(\mathbf R)$ is exactly the irreducible unitary representation $T_{\Omega_+}$ of $G=\Aff_0(\mathbf R)$ associated following the orbit method construction, to the orbit $\Omega_+$, which is the upper half-plane $\mathbf H \cong \mathbf R \rtimes \mathbf R^*$, i. e.
+$$(\exp(\hat{\ell}_Z)f)(y) = (T_{\Omega_+}(g)f)(y) = e^{iby}f(ay),\forall f\in L^2(\mathbf R^*, \frac{dy}{\vert y\vert}), $$ where  $g = \exp Z = \left(\begin{array}{cc} a & b \\ 0 & 1 \end{array}\right).$
\end{theorem}

By analogy, we have also
\begin{theorem}
The representation $\exp(\hat{\ell}_Z)$ of the group $G=\Aff_0(\mathbf R)$ is exactly the irreducible unitary representation $T_{\Omega_-}$ of $G=\Aff_0(\mathbf R)$ associated following the orbit method construction, to the orbit $\Omega_-$, which is the lower half-plane $\mathbf H \cong \mathbf R \rtimes \mathbf R^*$, i. e.
$$(\exp(\hat{\ell}_Z)f)(y) = (T_{\Omega_-}(g)f)(y) = e^{iby}f(ay),\forall f\in L^2(\mathbf R^*, \frac{dy}{\vert y\vert}), $$ where  $g = \exp Z = \left(\begin{array}{cc} a & b \\ 0 & 1 \end{array}\right).$
\end{theorem}
\begin{remark}{\rm
1. We have demonstrated how all the irreducible unitary representation of the connected group of affine transformations could be obtained from deformation quantization. It is reasonable to refer to the algebras of functions on co-adjoint orbits with this $\star$-product as {\it quantum ones}.

2. In a forthcoming work, we shall do the same calculation for the group of affine transformations of the complex straight line $\mathbf C$. This achieves the description of {\it quantum $\overline{MD}$ co-adjoint orbits}, see \cite{dndiep} for definition of $\overline{MD}$ Lie algebras.
}\end{remark}

\subsubsection{The group of affine transformations of the complex straight line}
Recall that the Lie algebra ${\mathfrak g} = \aff({\bf C})$ of affine transformations of the complex straight line is described as follows, see [D].

It is well-known that the group $\Aff({\bf C})$ is a four (real) dimensional Lie group which is isomorphism to the group of matrices:
$$\Aff({\bf C}) \cong  \left\{\left(\begin{array}{cc} a & b \\ 0 & 1 \end{array}\right) \vert a,b \in {\bf C}, a \ne 0 \right\}$$

The most easy method is to consider $X$,$Y$ as complex generators,
$X=X_1+iX_2$ and $Y=Y_1+iY_2$. Then from the relation $[X,Y]=Y$, we get$ [X_1,Y_1]-[X_2,Y_2]+i([X_1Y_2]+[X_2,Y_1]) = Y_1+iY_2$. 
This mean that the Lie algebra $\aff({\bf C})$ is a real 4-dimensional Lie algebra, having 4 generators with the only nonzero Lie brackets: $[X_1,Y_1] - [X_2,Y_2]=Y_1$; $[X_2,Y_1] + [X_1,Y_2] = Y_2$ and we can choose another basic noted again by the same letters to have more clear Lie brackets of this Lie algebra:
$$[X_1,Y_1] = Y_1; [X_1,Y_2] = Y_2; [X_2,Y_1] = Y_2; [X_2,Y_2] = -Y_1$$

\begin{remark}{\rm
The exponential map $$\exp: {\bf C}  \longrightarrow  {\bf  C}^{*} := {\bf C} \backslash \{0\}$$  giving by $z \mapsto e^z$ is just the covering map and therefore the universal covering of ${\bf} C^*$ is $\widetilde {\bf C}^* \cong {\bf C}$. As a consequence one deduces that $$\widetilde {\Aff}({\bf C}) \cong {\bf C} \ltimes{\bf C} \cong  \{(z,w) \vert z,w \in {\bf C} \}$$  with the following multiplication law:
$$(z,w)(z^{'},w^{'}) := (z+z',w+e^{z}w')$$
}\end{remark}

\begin{remark} {\rm 
The co-adjoint orbit of $\widetilde\Aff({\bf C})$ in ${\mathfrak g}^*$  passing through $F \in {\mathfrak g}^*$ is denoted by 
$$\Omega_{F} := K(\widetilde {\Aff}({\bf C})) F = \{K(g)F \vert g \in \widetilde \Aff({\bf C})\}$$
Then, (see [D]):
\begin{enumerate}
\item Each point $(\alpha,0,0,\delta)$ is 0-dimensional co-adjoint orbit $\Omega_{(\alpha,0,0,\delta)}$
\item The open set $\beta^{2}+\gamma^{2} \ne $ 0 is the single 4-dimensional co-adjoint orbit $\Omega_{F} = \Omega_{\beta^{2}+\gamma^{2} \ne 0} $. We shall also use $\Omega_{F}$ in form $\Omega_{F} \cong {\bf C} \times {{\bf C}}^*$.
\end{enumerate}
}\end{remark}

\begin{remark}{\rm 
Let us denote:
$$\mathbf H_{k} = \{w=q_{1}+iq_{2} \in {\bf C} \vert -\infty< q_1<+\infty ; 2k\pi < q_2< 2k\pi+2\pi\}; k=0,\pm1,\dots$$
$$L=\{{\rho}e^{i\varphi} \in {\bf C} \vert 0< \rho < +\infty; \varphi = 0\} \mbox{ and } {\bf C} _{k } = {\bf C} \backslash L$$
is a univalent sheet of the Riemann surface of the complex variable multi-valued analytic function $\Ln(w)$, ($k=0,\pm 1,\dots$)
Then there is a natural diffeomorphism $w \in \mathbf H_{k} \longmapsto e^{w} \in {\bf C}_k$ with each $k=0,\pm1,\dots.$ Now consider the map:
$${\bf C} \times {\bf C} \longrightarrow \Omega_F = {\bf C} \times {\bf C}^*$$
$$(z,w) \longmapsto (z,e^w),$$
with a fixed $k \in \mathbf Z$. We have a local diffeomorphism 
$$\varphi_k: {\bf C} \times {\bf H}_k \longrightarrow {\bf C} \times {\bf C}_k$$
 $$(z,w) \longmapsto (z,e^w) $$
This diffeomorphism $\varphi_k$ will be needed in the all sequel.
}\end{remark}

On ${\bf C}\times {\bf H}_k$ we have the natural symplectic form 
\begin{equation}\omega = \frac{1 }{2}[dz \wedge dw+d\overline {z} \wedge d\overline {w}],\end{equation} induced from $\mathbf C^2$.
Put $z=p_1+ip_2,w=q_1+iq_2$ and $(x^1,x^2,x^3,x^4)=(p_1,q_1,p_2,q_2) \in {\bf R}^4$, then
$$\omega = dp_1 \wedge dq_1-dp_2 \wedge dq_2.$$ The corresponding symplectic matrix of $\omega$ is 
$$ \wedge = \left(\begin{array}{cccc} 0 & -1 & 0 & 0 \\
		           1 & 0 & 0 & 0 \\
                                                0 & 0 & 0 & 1 \\
                                                0 & 0 & -1 & 0 \end{array}\right)
\mbox{   and   }          
 \wedge^{-1} = \left(\begin{array}{cccc} 0 & 1 & 0& 0 \\
		           -1 & 0 & 0 & 0 \\
                                                0 & 0 & 0 & -1 \\
                                                0 & 0 & 1 & 0 \end{array}\right)$$

We have therefore the Poisson brackets of functions as follows. With each $f,g \in {\bf C}^{\infty}(\Omega)$ 
$$\{f,g\} = \wedge^{ij}\frac{\partial f }{\partial x^i}\frac{\partial g}{\partial x^j} =  
\wedge^{12}\frac{\partial f }{\partial p_1}\frac{\partial g}{\partial q_1}+
\wedge^{21}\frac{\partial f }{\partial q_1}\frac{\partial g}{\partial p_1} +
\wedge^{34}\frac{\partial f}{\partial p_2}\frac{\partial g}{\partial q_2} +
\wedge^{43}\frac{\partial f}{\partial q_2}\frac{\partial g}{\partial p_2} = $$
$$\ \ \ \ \ \ \ =\frac{\partial f }{\partial p_1}\frac{\partial g}{\partial q_1} -
\frac{\partial f}{\partial q_1}\frac{\partial g}{\partial p_1} -
\frac{\partial f}{\partial p_2}\frac{\partial g }{\partial q_2} +
\frac{\partial f}{\partial q_2}\frac{\partial g }{\partial p_2} = $$
$$\ \ =2\Bigl[\frac{\partial f}{\partial z}\frac{\partial g}{\partial w} -
\frac{\partial f}{\partial w}\frac{\partial g}{\partial z} +
\frac{\partial f}{\partial \overline z}\frac{\partial g}{\partial \overline{w}} -
\frac{\partial f}{\partial \overline w}\frac{\partial g}{\partial \overline z}\Bigr]$$  

\begin{proposition}
Fixing the local  diffeomorphism $\varphi_k (k \in {\bf Z})$, we have: 
\begin{enumerate}
\item For any element $A \in \aff(\mathbf C)$, the corresponding Hamiltonian function $\widetilde{A}$  in local coordinates $(z,w)$ of the orbit $\Omega_F$  is of the form
$$\widetilde A\circ\varphi_k(z,w) = \frac{1}{2} [\alpha z +\beta e^w + \overline{\alpha} \overline{z} + \overline{\beta}e^{\overline {w}}]$$
\item In local coordinates $(z,w)$ of the orbit $\Omega_F$, the symplectic Kirillov form $\omega_F$ is just the standard form (1).
\end{enumerate}
\end{proposition}

{\it Computation of Operators $\hat{\ell}_A^{(k)}$.}

\begin{proposition}\label{Proposition 3.1}
With $A,B \in \aff({\bf C})$, the Moyal $\star$-product satisfies the relation:
\begin{equation} i \widetilde{A} \star i \widetilde{B} - i \widetilde{B} \star i \widetilde{A} = i[\widetilde{A,B} ]\end{equation}
\end{proposition}

For each $A \in \hbox{aff}({\bf C}$), the corresponding Hamiltonian function is 
$$\widetilde{A} = \frac{1}{2} [\alpha  z + \beta e^{w} + \overline{\alpha}\overline z + \overline{\beta} e^ {\overline w}] $$ 
and we can consider the operator  ${\ell}_A^{(k)}$  acting on dense subspace
$L^2({\bf R}^2\times ({\bf R}^2)^*,\frac{dp_1dq_1dp_2dq_2}{(2\pi)^2} )^{\infty}$
 of smooth functions by left $\star$-multiplication by $i \widetilde{A}$, i.e:
${\ell}_A^{(k)} (f) = i \widetilde{A} \star f$. Because of the relation in Proposition 3.1, we have
\begin{corollary}\label{Consequence 3.2}
 \begin{equation}{\ell}_{[A,B]}^{(k)} = {\ell}_A^{(k)} \star {\ell}_B^{(k)} - {\ell}_B^{(k)} \star {\ell}_A^{(k)} := {\Bigl[ {\ell}_A^{(k)},  {\ell}_B^{(k)}\Bigr]}^{\star}\end{equation}
\end{corollary}

From this it is easy to see that, the correspondence $A \in \aff({\bf C}) \longmapsto {\ell}_A^{(k)} = $i$\widetilde {A} \star$. is a representation of the Lie algebra $\aff({\bf C}$) on the space N$\bigl[[\frac{i}{2}]\bigr] $ of formal power series, see [G] for more detail.

Now, let us denote  ${\mathcal F}_z$(f) the partial Fourier transform of the function f from the variable $z=p_1+ip_2$ to the variable $\xi=\xi_1+i\xi_2$, i.e:
$${\mathcal F}_z(f )(\xi,w) = \frac{1 }{2\pi} \iint_{R^2} e^{-iRe(\xi \overline{z})} f(z,w)dp_1dp_2$$

Let us denote by $$\mathcal F_z^{-1}(f )(\xi,w) = \frac{1}{2\pi} \iint_{R^2} e^{iRe(\xi \overline{z})} f( \xi,w)d \xi_1d\xi_2$$ the inverse Fourier transform.

\begin{lemma}\label{Lemma 3.2}     Putting $g = g(z,w) = {\mathcal F}_z^{-1} (f )(z,w)$ we have:
\begin{enumerate}
\item
$$\partial_z g = \frac{i}{2}\overline\xi g \    ;  \partial_z^{r} g = {(\frac{i}{2}\overline\xi)}^r g, r=2,3,\dots $$
\item
 $$ \partial_{\overline z} g = \frac{i}{2}\xi g \    ;  \partial_{\overline z}^{r} g = {(\frac{i}{2}\xi)}^r g, r=2,3,\dots $$
\item 
$${\mathcal F}_z(zg) = 2i\partial_{\overline\xi}{\mathcal F}_z(g) = 2i\partial_{\overline \xi}f \    ;  {\mathcal F}_z(\overline{z}g) = 2i\partial_{\xi}{\mathcal F}_z(g) = 2i\partial_{\xi}f $$    
\item
$$ \partial_w g = \partial_w ({\mathcal F}_z^{-1}(f)) = {{\mathcal F}_z} ^{-1}(\partial_{w}f);\     \partial_{\overline w}g = \partial_{\overline w} ({\mathcal F}_z^{-1}(f)  = {{\mathcal F}_z}^{-1}(\partial_{\overline w}f)$$
\end{enumerate}  
\end{lemma}

We also need another Lemma which will be used in the sequel.

\begin{lemma}\label{Lemma 3.3}     With $g = {\mathcal F}_z^{-1}$$(f)($$z,w)$, we have: 
\begin{enumerate}
\item
$$    {\mathcal F}_z(P^0(\widetilde{A},g)) = i(\alpha \partial_{\overline \xi} + \overline {\alpha} \partial_{\xi})f + \frac{1}{ 2} \beta e^w f + \frac{1}{2} \overline {\beta} e^{\overline w} f. $$
\item
$$   {\mathcal F}_z(P^1(\widetilde{A},g)) = \overline {\alpha} \partial_{\overline w}f + \alpha \partial_{w}f - \overline {\beta} e^{\overline w} (\frac{i}{2} \xi)f - \beta e^w (\frac{i}{2}\overline {\xi})f. $$
\item
$$    {\mathcal F}_z(P^r(\widetilde{A},g)) = {(-1)}^r.2^{r-1}[\overline \beta{ e^{\overline w}} (\frac{i}{2}\xi)^r + \beta e^w (\frac{i}{2}\overline \xi)^r]f \ \ \ \ \ \      \forall r \ge 2. $$
\end{enumerate}
\end{lemma}

\begin{proposition}\label{Proposition 3.4}
For each $A = \left(\begin{matrix}\alpha & \beta \cr 0 & 0 \cr\end{matrix}\right) \in \aff({\bf C}) $
 and for each compactly supported $C^{\infty}$-function $f \in C_0^{\infty}({\bf C} \times {\bf H}_k)$, we have:
\begin{equation} {\ell}_A^{(k)}{f} := {\mathcal F}_z \circ \ell_A^{(k)} \circ {\mathcal F}_z^{-1}(f) = [\alpha (\frac{1}{2} \partial_w - \partial_{\overline \xi})f + \overline \alpha (\frac{1 }{2}\partial_{\overline w} - \partial_\xi)f + \end{equation}
$$+\frac{i}{2}(\beta e^{w-\frac{1}{2}\overline \xi} + \overline \beta e^{\overline w - \frac{1}{2} \xi})f] $$
\end{proposition}

\begin{remark}\label{Remark 3.5}{\rm Setting new variables  u = $w - \frac{1}{ 2}\overline{\xi}$;$v = w + \frac{1 }{2}{\overline{\xi}}$ we have
\begin{equation}\hat {\ell}_A^{(k)}(f) = \alpha\frac{ \partial f }{\partial u}+ \overline{\alpha}\frac{\partial f }{\partial{\overline{u}}}+ \frac{i }{2}(\beta e^{u}+\overline{\beta}e^{\overline{u}})f \vert_{(u,v)}\end{equation}
i.e $\hat {\ell}_A^{(k)} = \alpha\frac{ \partial }{\partial u}+ \overline{\alpha}\frac{ \partial  }{\partial{\overline{u}}}+ \frac{i }{2}(\beta e^{u}+\overline{\beta}e^{\overline{u}})$,which provides a ( local) representation of the Lie algebra  aff({\bf C}).
}\end{remark}

{\it The Irreducible Representations of $\widetilde{\Aff}({\bf C})$. }
Since $\hat {\ell}_A^{(k)}$ is a representation of the Lie algebra  $\widetilde{\hbox {Aff}} ({\bf C})$, we have:
$$\exp(\hat {\ell}_A^{(k)}) = \exp\bigl(\alpha\frac{ \partial }{\partial {u}}+ \overline{\alpha}\frac{ \partial  }{\partial{\overline{u}}}+ \frac{i }{2}(\beta e^{u}+\overline{\beta}e^{\overline{u}})\bigr)$$ is just the corresponding representation of the corresponding connected and simply connected Lie group $\widetilde\Aff ({\bf C})$.

Let us first recall the well-known list of all the irreducible unitary representations of the group of affine transformation of the complex straight line, see [D] for more details.

\begin{theorem}\label{Theorem 4.1}
Up to unitary equivalence, every irreducible unitary representation of $\widetilde{\hbox {Aff}} ({\bf C})$ is unitarily equivalent to one the following one-to-another nonequivalent irreducible unitary representations:
\begin{enumerate}
\item The unitary characters of the group, i.e the one dimensional unitary representation $U_{\lambda},\lambda \in {\bf C}$, acting in ${\bf C}$ following the formula
$U_{\lambda}(z,w) = e^{{i\Re(z\overline{\lambda})}}, \forall (z,w) \in \widetilde{\Aff} ({\bf C}), \lambda \in {\bf C}.$
\item The infinite dimensional irreducible representations $T_{\theta},\theta \in {\mathbf S}^1$, acting on the Hilbert space $L^{2}(\mathbf R\times \mathbf S ^1)$ following the formula:
\begin{equation}\Bigr[T_{\theta}(z,w)f\Bigl](x) = \exp \Bigr(i(\Re(wx)+2\pi\theta[\frac{\Im(x+z) }{2\pi}])\Bigl)f(x\oplus z),\end{equation}
Where \ $(z,w) \in\widetilde{\Aff}({\bf C})$  ;  $x \in {\bf R}\times {\mathbf S} ^1= {\bf C} \backslash \{0\}; f \in L^{2}({\bf R}\times {\mathbf S} ^1);$
$$ x\oplus z = Re(x+z) +2 \pi i \{\frac{\Im(x+z) }{2\pi}\}$$
\end{enumerate}
\end{theorem}
In this section we will prove the following important Theorem which
is very interesting for us both in theory and practice.
\par
\begin{theorem}\label{Theorem 4.2}
The representation $\exp(\hat {\ell}_A^{(k)})$ of the group $\widetilde{\Aff}({\bf C})$ is the irreducible unitary representation 
$T_\theta$ of $\widetilde{\Aff}({\bf C})$ associated, following the orbit method construction, to the orbit $\Omega$, i.e:
$$\exp(\hat {\ell}_A^{(k)})f(x) = [T_\theta (\exp A)f](x),$$
where $f \in L^{2}({\bf R}\times {\mathbf S} ^1) ; A = \begin{pmatrix}\alpha & \beta \cr 0 & 0 \cr\end{pmatrix} \in \aff({\bf C}) ; \theta \in {\mathbf S}^1 ; k = 0, \pm1,\dots$
\end{theorem} 

\begin{remark}\label{Remark 4.3} {\rm
We say that a real Lie algebra ${\mathfrak g}$ is in the class $\overline{MD}$ if every K-orbit is of dimension, equal 0 or dim ${\mathfrak g}$. Further more, one proved that
([D, Theorem 4.4]) 
Up to isomorphism, every Lie algebra of class $\overline {MD}$ is one of the following:
\begin{enumerate}
\item Commutative Lie algebras.
\item Lie algebra $\aff({\bf R})$ of affine transformations of the real straight line
\item Lie algebra $\aff({\bf C})$ of affine transformations of the complex straight line.
\end{enumerate}
Thus, by calculation for the group of affine transformations of the real straight line $\Aff({\bf R})$ in [DH] and here for the group affine transformations of the complex straight line $\Aff({\bf C})$ we obtained  a description of the quantum $\overline {MD}$ co-adjoint orbits.
}\end{remark}

\subsection{$MD_4$-groups}
We refer the reader to the results of Nguyen Viet Hai \cite{hai3}-\cite{hai4} for the class of $MD_4$-groups (i.e. 4-dimensional solvable Lie groups, all the coadjoint of which are of dimension 0 or maximal). It is interesting that here he obtained the same exact computation for $\star$-products and all representations.

\subsection{$SO(3)$} As an typical example of compact Lie group, the author proosed Job A. Nable to consider the case of $SO(3)$. We refer the reader to the  results of Job Nable \cite{nable1}-\cite{nable3}. In these examples, it is interesting that the $\star$-products, in some how as explained in these papers, involved the Maslov indices and Monodromy Theorem.

\section*{Acknowledgments} 

This work was completed during the stay of the first author as a visiting
mathematician at the Department of mathematics, The University of Iowa. The author would like to express the deep and sincere thanks to Professor Tuong Ton-That and his spouse, Dr. Thai-Binh Ton-That for their effective helps and kind attention they provided during the stay in Iowa, and also for a discussion about the PBW Theorem. The deep thanks are also addressed to the organizers of the Seminar on Mathematical Physics, Seminar on Operator Theory in Iowa and the Iowa-Nebraska Functional Analysis Seminar (INFAS), in particular the professors Raul Curto, Palle Jorgensen, Paul Muhly and Tuong Ton-That for the stimulating scientific atmosphere.
 
The author would like to thank the University of Iowa for the hospitality and the scientific support, the Alexander von Humboldt Foundation, Germany, for an effective support.

\end{document}